\documentclass[12pt,reqno,a4paper]{amsart}
\usepackage{amssymb,amscd, hyperref, epsf}

\begin{document}

\let\kappa=\varkappa
\let\epsilon=\varepsilon
\let\phi=\varphi
\let\p\partial

\def\Z{\mathbb Z}
\def\R{\mathbb R}
\def\C{\mathbb C}
\def\Q{\mathbb Q}
\def\P{\mathbb P}
\def\HH{\mathrm{H}}
\def\ss{X}

\def\conj{\overline}
\def\Beta{\mathrm{B}}
\def\const{\mathrm{const}}
\def\ov{\overline}
\def\wt{\widetilde}
\def\wh{\widehat}

\renewcommand{\Im}{\mathop{\mathrm{Im}}\nolimits}
\renewcommand{\Re}{\mathop{\mathrm{Re}}\nolimits}
\newcommand{\codim}{\mathop{\mathrm{codim}}\nolimits}
\newcommand{\id}{\mathop{\mathrm{id}}\nolimits}
\newcommand{\Aut}{\mathop{\mathrm{Aut}}\nolimits}
\newcommand{\lk}{\mathop{\mathrm{lk}}\nolimits}
\newcommand{\sign}{\mathop{\mathrm{sign}}\nolimits}
\newcommand{\rk}{\mathop{\mathrm{rk}}\nolimits}
\def\Jet{{\mathcal J}}
\def\FC{{\mathrm{FCrit}}}
\def\sS{{\mathcal S}}
\def\lcan{\lambda_{\mathrm{can}}}
\def\ocan{\omega_{\mathrm{can}}}

\renewcommand{\mod}{\mathrel{\mathrm{mod}}}

\newtheorem{mainthm}{Theorem}
\newtheorem{thm}{Theorem}[subsection]
\newtheorem{lem}[thm]{Lemma}
\newtheorem{prop}[thm]{Proposition}
\newtheorem{cor}[thm]{Corollary}

\theoremstyle{definition}
\newtheorem{exm}[thm]{Example}
\newtheorem{rem}[thm]{Remark}
\newtheorem{df}[thm]{Definition}

\renewcommand{\thesubsection}{\arabic{subsection}}
\numberwithin{equation}{subsection}

\title{Non-negative Legendrian isotopy in $ST^*M$}
\author[Chernov \& Nemirovski]{Vladimir Chernov and Stefan Nemirovski}
\thanks{The second author was supported by grants from DFG, RFBR, Russian Science Support Foundation,
and the programme ``Leading Scientific Schools of Russia.''}
\address{Department of Mathematics, 6188 Kemeny Hall,
Dartmouth College, Hanover, NH 03755, USA}
\email{Vladimir.Chernov@dartmouth.edu}
\address{%
Steklov Mathematical Institute, 119991 Moscow, Russia;\hfill\break
\strut\hspace{8 true pt} Mathematisches Institut, Ruhr-Universit\"at Bochum, 44780 Bochum, Germany}
\email{stefan@mi.ras.ru}

\begin{abstract}
It is shown that if the universal cover of a manifold $M$ is an open manifold,
then two different fibres of the spherical cotangent bundle $ST^*M$
cannot be connected by a non-negative Legendrian isotopy.
This result is applied to the study of causality in globally hyperbolic spacetimes.
It is also used to strengthen a result of Eliashberg, Kim, and Polterovich
on the existence of a partial order on $\widetilde{\mathrm{Cont}}_0 (ST^*M)$.
\end{abstract}

\maketitle

\subsection{Introduction}
Let $M$ be a connected not necessarily orientable manifold of dimension $m\ge 2$
and let $\pi_M:ST^*M\to M$ be its spherical cotangent bundle.
It is well-known that $ST^*M$ carries a canonical co-oriented contact structure.
An isotopy $\{L_t\}_{t\in[0,1]}$ of Legendrian submanifolds in a co-oriented
contact manifold is called {\it non-negative\/} if it can be parameterised
in such a way that the tangent vectors of the trajectories of individual points
lie in the non-negative tangent half-spaces defined by the contact structure,
see Definition~\ref{defnonneg}. For a generic Legendrian isotopy in $ST^*M$,
this property can be expressed by saying that the co-oriented wave fronts
$\pi_M(L_t)\subset M$ move in the direction of their co-orientation.

\begin{thm}
\label{mainresult}
Assume that the universal cover of $M$ is an open manifold.
Then there does not exist a non-negative Legendrian isotopy
connecting two different {\rm (}nonoriented\/{\rm )} fibres of $ST^*M$.
\end{thm}

In the special case when $M$ can be covered by on open subset of~$\R^m$,
this statement was proved by Colin, Ferrand, and Pushkar~\cite{CFP}.
Independently, a slightly stronger result was obtained by the present authors in the course
of the proof of the so-called Legendrian Low Conjecture from Lorentz geometry~\cite{ChNe}.
Theorem~\ref{mainresult} allows us to extend the results
of~\cite{ChNe} to a wider class of Lorentz manifolds, see~\S\ref{lowapplic}.

A closely related notion of non-negative contact isotopy plays a key role
in the orderability problem for contactomorphism groups, see~\cite{ElPo},
\cite{ElKiPo}, and~\cite{Bh}. Theorem~\ref{mainresult} can be applied
to settle a question left open in~\cite{ElKiPo}, see~\S\ref{orderapplic}.

It is easy to see that the assertion of Theorem~\ref{mainresult} is false if $M$ carries
a Riemann metric turning it into a $Y_\ell^x$-manifold, see Example~\ref{Yxl}.
In particular, it is false if $M$ is a metric quotient of the standard sphere.
Hence, the hypothesis of Theorem~\ref{mainresult} cannot be weakened for surfaces
(this is obvious) and $3$-manifolds (this follows from Perelman's work
on the Poincar\'e conjecture~\cite{Perelman1,Perelman2,Perelman3}).
On the other hand, it seems very likely
that the result can be strengthened in all higher dimensions.

The proof of Theorem~\ref{mainresult} is based on Viterbo's invariants
of generating functions~\cite{Vi}. However, it is different from
the arguments in~\cite{CFP} and~\cite{ChNe} already in the case when~$M=\R^m$.
In particular, no use is made of the identification $ST^*\R^m\cong\Jet^1(S^{m-1})$.

All manifolds, maps etc.~are assumed to be smooth unless the opposite is explicitly stated,
and the word {\it smooth\/} means $C^{\infty}$.

\subsection{Non-negative Legendrian isotopies}
\label{legisotop}
Let $(Y,\ker\alpha)$ be a contact manifold with a co-oriented
contact structure defined by a contact form~$\alpha$.

\begin{df}
\label{defnonneg}
A Legendrian isotopy $\{L_t\}_{t\in[0,1]}$ in $(Y,\ker\alpha)$ is called {\it non-negative\/}
if it has a parameterisation $F:L_0\times[0,1]\to Y$ such that
$(F^*\alpha)\left(\frac{\p}{\p t}\right)\ge 0$.
If the latter inequality is strict, the isotopy is called {\it positive}.
\end{df}

Clearly, this definition does not depend on the choice of the parameterisation $F$
of the Legendrian isotopy and on the choice of the contact form defining the
co-oriented contact structure. It is also obvious that (co-orientation preserving)
contactomorphisms preserve the property of being non-negative or positive.

\begin{lem}
\label{posapprox}
Let $\{L_t\}_{t\in[0,1]}$ be a non-negative Legendrian isotopy of compact submanifolds
such that $L_0\cap L_1=\varnothing$. Then there exists a {\rm (}$C^\infty$-close\/{\rm )}
positive Legendrian isotopy with the same ends.
\end{lem}

\begin{proof}[Sketch of proof]
Let $L_t'=\psi_{\epsilon t}(L_t)$, where $\psi_t$ is the Reeb flow on $(Y,\alpha)$
and $\epsilon>0$. The isotopy $\{L_t'\}$ is positive. If $\epsilon$ is small enough,
then $L_1'$ is Legendrian isotopic to $L_1$ in $Y\setminus L_0$ and therefore there
exists a contactomorphism $\phi$ supported in $Y\setminus L_0$ such that $\phi(L_1')=L_1$.
Thus, $L''_t:=\phi(L_t')$ is a positive Legendrian isotopy connecting $L_0$ and $L_1$.
\end{proof}

An advantage of positive isotopies is that positivity is an open
condition and hence one can make a positive Legendrian isotopy generic by a
small perturbation.

\begin{df}
\label{genpos}
An isotopy $\{L_t\}_{t\in[0,1]}$ is in general position with respect
to a submanifold $\Lambda\subset Y$ of codimension $\dim_\R L_0$ if
it has a parameterisation $F:L_0\times [0,1]\to Y$ such that
\begin{itemize}
\item[a)] $F^{-1}(\Lambda)$ is a $1$-dimensional submanifold in $L_0\times [0,1]$;
\item[b)] the projection $F^{-1}(\Lambda)\to [0,1]$ has isolated critical points;
\item[c)] $F^{-1}(\Lambda)$ is transverse to $L_0\times\{0\}$ and $L_0\times\{1\}$.
\end{itemize}
\end{df}

Note that a point $(x,\tau)\in F^{-1}(\Lambda)$ which is not critical
for the projection $F^{-1}(\Lambda)\to [0,1]$ lies on the graph of a section
of this projection over a non-trivial closed interval $[t',t'']\ni\tau$.
In other words, there exists a curve $\gamma:[t',t'']\to L_0$ such that
$\gamma(\tau)=x$ and $F(\gamma(t),t)\in \Lambda$ for all $t\in[t',t'']$.

\subsection{Exact pre-Lagrangian submanifolds}
Suppose that $\Lambda$ is an $m$-dimensional submanifold
of a $(2m-1)$-dimensional contact manifold $(Y,\ker\alpha)$
such that
\begin{equation}
\label{exprelag}
df=e^h\alpha|_\Lambda
\end{equation}
for some functions $f,h:\Lambda\to\R$.
Then $\Lambda$ is said to be {\it exact pre-Lagrangian\/} and
$f$ is called a {\it contact potential\/} on~$\Lambda$. (The
terminology will be explained in \S\ref{sympl}.)

The following lemma shows that the contact potential is non-decreasing
along any curve traced on $\Lambda$ by a non-negative Legendrian isotopy.

\begin{lem}
\label{fmonoton}
Let $L_t=\phi_t(L_0)$, $t\in[0,1]$, be a non-negative Legendrian isotopy.
Suppose that $\gamma:[0,1]\to L_0$ is a curve such that $\phi_t(\gamma(t))\in\Lambda$,
where $\Lambda$ is an exact pre-Lagrangian submanifold with contact potential~$f$.
Then the function $t\mapsto f(\phi_t(\gamma(t)))$ is non-decreasing.
\end{lem}

\begin{proof}
This follows from the definitions and the chain rule. Indeed,
\begin{eqnarray*}
\frac{d}{dt}f(\phi_t(\gamma(t)))&=& df\bigl(\frac{d\phi_t}{dt}(\gamma(t))+d\phi_t(\gamma(t))\,\frac{d\gamma}{dt}(t)\bigr)\\[2pt]
&=&e^h\left[\alpha\bigl(\frac{d\phi_t}{dt}(\gamma(t)\bigr)+\alpha\bigl(d\phi_t(\gamma(t))\,\frac{d\gamma}{dt}(t)\bigr)\right].
\end{eqnarray*}
The first summand in square brackets is non-negative by the definition of non-negative
Legendrian isotopy and the second one is zero because the isotopy is Legendrian
and $\frac{d\gamma}{dt}$ is tangent to~$L_0$. Hence, the derivative of our function is non-negative.
\end{proof}

\subsection{Symplectisation}
\label{sympl}
Let $(Y,\ker\alpha)$ be a contact manifold. Its symplectisation $Y^{\mathrm{symp}}$
is the (exact) symplectic manifold $(Y\times\R,d(e^s\alpha))$.

\begin{exm}
\label{sympspher}
Let $Y\subset T^*M$ be the unit sphere bundle with respect to
a Riemann metric on~$M$. Then $\alpha:={\lcan|}_Y$ is a contact
form defining the canonical contact structure on $Y\cong ST^*M$. The map
$$
Y^{\mathrm{symp}} \ni (\xi,s)\longmapsto e^s\xi\in T^*M
$$
is a symplectomorphism onto the complement of the zero section of $T^*M$
such that the pull-back of the canonical $1$-form $\lcan$ is precisely $e^s\alpha$.\qed
\end{exm}

A contactomorphism $\phi:Y\to Y$ lifts to a symplectomorphism
$\wt\phi:Y^{\mathrm{symp}}\to Y^{\mathrm{symp}}$ defined by the formula
$$
\wt\phi(x,s):=(\phi(x),s-\eta(x)),
$$
where $\eta:Y\to\R$ is the function such that $\phi^*\alpha=e^\eta\alpha$.
It follows from this definition that if $\{\phi_t\}_{t\in[0,1]}$ is a contact isotopy of $Y$,
then $\{\wt\phi_t\}_{t\in[0,1]}$ is a {\it Hamiltonian\/} isotopy of $Y^{\mathrm{symp}}$.
(It can be defined by the Hamiltonian function $\wt H(x,s,t)=-\alpha(\frac{d\phi_t}{dt})$.)

An exact pre-Lagrangian submanifold $\Lambda$ with a contact potential $f$
such that $df=e^h\alpha|_\Lambda$ lifts to an exact Lagrangian submanifold
$$
\wt\Lambda=\{(x,h(x))\in Y^{\mathrm{symp}}\mid x\in\Lambda\}\subset Y^{\mathrm{symp}}.
$$
Indeed, the function $\wt f:\wt\Lambda\to\R$, $\wt f(x,h(x))=f(x)$,
is a primitive for the $1$-form $e^s\alpha|_{\wt\Lambda}$.
Note that for any contactomorphism $\phi:Y\to Y$, the image $\phi(\Lambda)$
with the contact potential $f\circ \phi^{-1}$ is exact pre-Lagrangian
and $\wt{\phi(\Lambda)}=\wt\phi(\wt\Lambda)$.

\subsection{Generating functions}
Let $M$ be a manifold, which may be open or have boundary.
Consider the product $M\times\R^N$ for some $N\ge 0$ and let
$\pi:M\times\R^N\to M$ be the projection onto~$M$.
For a function $\sS:M\times\R^N\to\R$, consider the set of its
fibre critical points
$$
\FC(\sS):=\{z\in M\times\R^N \mid d\sS(z)|_{\R^N}=0\}.
$$
Note that there is a natural fibrewise map
$$
d_M \sS:\FC(\sS)\to T^*M
$$
which associates to a point $z\in \FC(\sS)$ the linear form
$v\mapsto d\sS(z)(\widehat v)$ on $T_{\pi(z)}M$, where $\widehat v\in T_z (M\times\R^N)$
is any tangent vector such that $d\pi(\widehat v)=v\in T_{\pi(z)}M$.

A function $\sS:M\times\R^N\to\R$ is called a {\it generating function}
for a Lagrangian submanifold $L\subset T^*M$ if it satisfies
the following two conditions:
\begin{itemize}
\item[(GF1)] its set of fibre critical points is cut out transversely;
\item[(GF2)] the map $d_M \sS:\FC(\sS)\to T^*M$ is a diffeomorphism onto~$L$.
\end{itemize}
Note that $\sS\circ (d_M \sS)^{-1}:L\to\R$ is a primitive for ${\lcan|}_L$.
Hence, a Lagrangian submanifold of $T^*M$ admitting a generating function is exact.

A generating function $\sS:M\times\R^N\to\R$ is called {\it quadratic at infinity\/} if furthermore
\begin{itemize}
\item[(GF3)] $\sS(y,\xi)=\sigma(y,\xi)+Q(\xi)$, where $Q$ is a non-degenerate quadratic
form on $\R^N$ and the projection $\pi:\mathop{\mathrm{supp}}\sigma\to M$ is a proper map.
\end{itemize}
Note that a Lagrangian submanifold $L\subset T^*M$ admitting a quadratic at infinity generating function is {\it properly embedded}, i.\,e., the projection $L\to M$ is a proper map.

\begin{prop}
\label{hamhom}
Let $\{L_t\}_{t\in[0,1]}$ be a compactly supported isotopy of
properly embedded Lagrangian submanifolds in $T^*M$.
Suppose that
\begin{itemize}
\item[a)] $L_0$ admits a quadratic at infinity generating function\/{\rm ;}
\item[b)] there exists a family of functions $f_t:L_t\to\R$ such that $df_t={\lcan|}_{L_t}$.
\end{itemize}
Then there exists a family $\sS_t:M\times\R^N\to\R$
of quadratic at infinity generating functions for~$L_t$
such that $\sS_t\circ (d_M \sS_t)^{-1}=f_t$ for all $t\in [0,1]$.
\end{prop}

\begin{proof}
This is a minor extension of the Laudenbach--Sikorav theorem~\cite{LaSi}.
It can be obtained, for instance, by applying the version of Chekanov's theorem~\cite{Ch}
for properly embedded Legendrian submanifolds~\cite[Sec.~4]{ElGr}
to the Legendrian isotopy
$$
\wh L_t:=\{(x,f_t(x))\in\Jet^1(M)\mid x\in L_t\}
$$
in the $1$-jet bundle of~$M$.
\end{proof}

\subsection{Critical values of quadratic at infinity functions}
Let $S:\R^N\to\R$ be a function quadratic at infinity in the sense that
$$
S(z)=\sigma(z)+Q(z),
$$
where $\sigma$ has compact support and $Q$ is a non-degenerate quadratic form on~$\R^N$.
(We will eventually take $S$ to be the restriction of a quadratic
at infinity generating function $\sS:M\times\R^N\to\R$ to the fibre $\{x\}\times\R^N$
over a point $x\in M$.)
Following Viterbo~\cite[\S2]{Vi}, let us define an invariant $c_-(S)\in\R$
of such a function.

Consider the sublevel sets
$$
S^c:=\{z\in\R^N\mid S(z)\le c\}
$$
and denote by $S^{-\infty}$ the set $S^c$ for a sufficiently negative $c\ll 0$.
Pick a $Q$-negative linear subspace $V\subset\R^N$ of maximal possible
dimension~$\kappa$. The relative homology class
$[V]\in\HH_\kappa(\R^N,S^{-\infty})$ does not depend on the choice of $V$.
Set
$$
c_-(S):=\inf\{c\in\R\mid [V]\in \imath_*\HH_\kappa(S^c,S^{-\infty})\},
$$
where $\imath_*:\HH_\kappa(S^c,S^{-\infty})\to \HH_\kappa(\R^N,S^{-\infty})$
is the homomorphism of relative homology groups
induced by the inclusion $\imath:S^c\to \R^N$.

By Morse theory, $c_-(S)$ is a critical value of $S$. In particular,
if $S$ has a single critical point $z_0\in\R^N$, then $c_-(S)=S(z_0)$.

We shall need the following version of Viterbo's monotonicity lemma~\cite[Lemma 4.7]{Vi}
adapted to our situation, cf.\ Lemma~\ref{fmonoton}.

\begin{lem}
\label{Smonoton}
Let $\{S_t\}_{t\in[0,1]}$ be a family of quadratic at infinity functions on~$\R^N$
and let
$$
C:=\{(z,\tau)\in \R^N\times [0,1]\mid dS_\tau(z)=0\}=\bigcup_{\tau\in [0,1]} \mathrm{Crit}(S_\tau)\times\{\tau\}.
$$
Suppose that for any $(z,\tau)$ from a dense subset $C'\subseteq C$ there exists
a non-trivial closed interval $[t',t'']\ni \tau$ and a curve $\gamma:[t',t'']\to \R^N$
such that
\begin{itemize}
\item[a)] $\gamma(\tau)=z;$
\item[b)] $\gamma(t)\in\mathrm{Crit}(S_t)$ for all $t\in[t',t''];$
\item[c)] the function $t\mapsto S_t(\gamma(t))$ is non-decreasing on~$[t',t'']$.
\end{itemize}
Then $t\mapsto c_-(S_t)$ is a non-decreasing {\rm (}continuous\/{\rm )} function on $[0,1]$.
\end{lem}

\begin{proof}
According to~\cite[Lemma 4.7]{Vi}, the claim will follow if we show that
$\frac{\p S_\tau}{\p t}(z)\ge 0$ for all $(z,\tau)\in C$.
If $(z,\tau)\in C'$, we have
$$
0\le \frac{d}{dt} S_t(\gamma(t))|_{t=\tau} = \frac{\p S_\tau}{\p t}(z)
+ dS_\tau(z)\, \frac{d\gamma}{dt}(\tau)= \frac{\p S_\tau}{\p t}(z),
$$
where $\gamma$ is a curve satisfying conditions (a)--(c) for $\tau$ and~$z$.
Since $C'$ is dense in $C$, this inequality is valid for all $(z,\tau)\in C$.
\end{proof}

\subsection{The Importance of Being Open}
\label{constructionM}
Let $M$ be an open manifold. We identify $ST^*M$ with the unit sphere
bundle in $T^*M$ for some Riemann metric on $M$ and view the complement
to the zero section of $T^*M$ as the symplectisation of $ST^*M$,
see Example~\ref{sympspher}. Let $\pi_M$ denote the bundle projection
$T^*M\to M$.

\begin{lem}
There exists a function $\Phi:M\to\R$ without critical points.
\end{lem}

\begin{proof}
This is well-known, see~\cite[Lemma 1.15]{Go}.
\end{proof}

\begin{df}
\label{maindef}
Let
$\Lambda_\Phi:= \left\{ \dfrac{d\Phi(x)}{\|d\Phi(x)\|} \mid x\in M\right\}\subset ST^*M.$
\end{df}

It is clear from the definition of the canonical $1$-form that
$\Lambda_\Phi$ is an exact pre-Lagrangian submanifold of $ST^*M$
and the function
$$
f_\Phi(\zeta)=\Phi(\pi_M(\zeta))
$$
is a contact potential on $\Lambda_\Phi$. The associated Lagrangian lift
$$
\wt\Lambda_\Phi=\left\{ d\Phi(x) \mid x\in M\right\}\subset T^*M
$$
is just the graph of the differential of $\Phi$. It has an obvious
generating function
$$
\sS_\Phi:M\times\R^0\to\R, \quad \sS_\Phi(x\times\{\mathrm{pt}\}):=\Phi(x).
$$

\subsection{Proof of Theorem~\ref{mainresult}}
Let $M$ be a manifold (universally) covered by an open manifold. Suppose that there
exists a non-negative Legendrian isotopy $\{L_t\}_{t\in[0,1]}$ connecting
two different fibres of $ST^*M$. Since such an isotopy lifts to the spherical
cotangent bundle of the covering manifold,  we may assume that $M$ is itself
an open manifold. By Lemma~\ref{posapprox}, we may also assume that the
Legendrian isotopy is positive.

Let $\Lambda_\Phi$ be the exact pre-Lagrangian submanifold
with contact potential $f_\Phi=\Phi\circ\pi_M$ defined in~\S\ref{constructionM}.
Applying a global contactomorphism induced by a suitable diffeomorphism of~$M$,
we can arrange that $L_0=ST_{x_0}^*M$ and $L_1=ST^*_{x_1} M$, where the points
$x_0,x_1\in M$ are such that $\Phi(x_0)>\Phi(x_1)$. Furthermore, we can put
the isotopy in general position with respect to $\Lambda_\Phi$ in the
sense of Definition~\ref{genpos}, leaving $L_0$ and $L_1$ fixed
(because they are already transversal to $\Lambda_\Phi$).

Let $\{\phi_t\}_{t\in[0,1]}$ be a compactly supported contact isotopy
of $ST^*M$ such that $L_t=\phi_t(L_0)$ for all~$t\in[0,1]$.
(Such an isotopy exists by the Legendrian isotopy extension theorem.)
Consider the Hamiltonian isotopy of exact Lagrangian submanifolds
$(\wt\phi_t)^{-1}(\wt\Lambda_\Phi)\subset T^*M$ and the functions
$\wt f_\Phi\circ\wt\phi_t$ on these manifolds, see \S\ref{sympl}.
By Proposition~\ref{hamhom}, there exists a family of quadratic at
infinity generating functions
$$
\sS_t:M\times\R^N\to\R
$$
for $(\wt\phi_t)^{-1}(\wt\Lambda_\Phi)\subset T^*M$ such that
$$
\sS_t\circ (d_M\sS_t)^{-1}=\wt f_\Phi\circ\wt \phi_t.
$$
Let
$$
S_t:=\sS_t(x_0,\cdot): \R^N\to \R
$$
be the restrictions of $\sS_t$ to the fibre $\{x_0\}\times\R^N$.
By construction, the map
$$
\mathop{\mathrm{Crit}}(S_t)\ni z
\longmapsto
\phi_t\left(\frac{d_M\sS_t(z)}{\|d_M\sS_t(z)\|}\right)\in L_t
$$
establishes a bijective correspondence between the set of critical
points of $S_t$ and the intersection $L_t\cap\Lambda_\Phi$.
Furthermore, the value of $S_t$ at a point $z\in \mathop{\mathrm{Crit}}S_t$
is equal to the value of $f_\Phi$ at the corresponding point in~$\Lambda_\Phi$.

In particular, $S_0$ and $S_1$ each have a single critical point corresponding
to the intersection of $\Lambda_\Phi$ with $L_0=ST^*_{x_0}M$ and
$L_1=ST^*_{x_1}M$, respectively. Since $f_\Phi=\Phi\circ\pi_M$,
we see that
$$
c_-(S_0)=\Phi(x_0)\quad\mbox{ and }\quad c_-(S_1)=\Phi(x_1).
$$
Hence,
\begin{equation}
\label{contra}
c_-(S_0)>c_-(S_1)
\end{equation}
by our choice of the points $x_0$ and $x_1$.

On the other hand, it follows from Lemma~\ref{fmonoton} and the discussion after
Definition~\ref{genpos} that the family of functions $\{S_t\}_{t\in[0,1]}$ satisfies the
hypotheses of Lemma~\ref{Smonoton}. Thus, $c_-(S_t)$ is a non-decreasing
function of $t$ and therefore $c_-(S_0)\le c_-(S_1)$, which contradicts~(\ref{contra}).

This contradiction shows that a non-negative Legendrian isotopy cannot connect
two different fibres of $ST^*M$.\qed

\begin{cor}
\label{loops}
If the universal cover of $M$ is an open manifold, then there does not exist
a positive Legendrian loop in the Legendrian isotopy class of the fibre of $ST^*M$.
\end{cor}

\begin{proof}
Suppose that $\{L_t\}_{t\in[0,1]}$ is a positive Legendrian isotopy such that
$L_0=L_1=ST^*_{x}M$. If $\{\phi_t\}_{t\in[0,1]}$ is a contact isotopy
such that $\phi_0=\mathrm{id}$, $\phi_1(ST^*_{x}M)=ST^*_{y}M$ for some $y\ne x$,
and $\|\frac{d\phi_t}{dt}\|$ is sufficiently small, then the isotopy $\{\phi_t(L_t)\}_{t\in[0,1]}$
is positive and connects two different fibres of $ST^*M$, which contradicts Theorem~\ref{mainresult}.
\end{proof}

\begin{rem}
\label{importantremark}
With a little more work, it can be shown that any {\it non-negative\/} Legendrian loop
in the Legendrian isotopy class of the fibre must be constant, cf.~\cite[Corollaries~5.5 and~6.2]{ChNe}.
\end{rem}

\begin{exm}
\label{Yxl}
Suppose that there exists a Riemann metric $g$ on $M$ such that $(M,g)$ is a $Y_\ell^x$-manifold
for some $x\in M$ and $\ell>0$, i.\,e., such that all $g$-geodesics starting from $x$ return to $x$
in time~$\ell$, see~\cite[Definition~7.7(c)]{Besse}. Then moving the fibre $ST^*_x M$ along the
(co-)geodesic flow on $ST^*M$ defines a positive Legendrian loop based at $ST^*_x M$.
Thus, Corollary~\ref{loops} and Theorem~\ref{mainresult} do not hold for such a manifold~$M$.

Note that if $\dim M=2$ or~$3$, then either the universal cover of $M$ is open
or $M$ admits a Riemann metric turning it into a $Y_\ell^x$-manifold.
For $\dim M=2$, this statement follows immediately from the classification of surfaces.
For $\dim M=3$, the Poincar\'e conjecture proved by Perelman~\cite{Perelman1,Perelman2,Perelman3}
implies that the universal cover of~$M$ is either non-compact or diffeomorphic to~$S^3$.
In the latter case, the elliptisation conjecture also proved by Perelman guarantees
that $M$ is diffeomorphic to a quotient of the standard round $S^3$
by the action of a finite group of isometries and the quotient metric turns $M$ into
a $Y_\ell^x$-manifold. Thus, Theorem~1.1 fails for every surface or $3$-manifold
such that its universal cover is not open.

The weak form of the Bott--Samelson theorem proved by B\'erard-Bergery, see~\cite{BerardBergery}
and~\cite[Theorem 7.37]{Besse}, says that if $(M,g)$ is a $Y_\ell^x$-manifold,
then the universal cover of $M$ is compact and the rational cohomology ring $H^*(M, \Q)$ is generated by one element.
In view of the preceding discussion, it seems natural to ask whether the latter property
is also shared by all manifolds $M$ such that there exists a positive Legendrian loop
in the Legendrian isotopy class of the fibre of $ST^*M$.\qed
\end{exm}

\begin{rem}\label{moregeneralbundles}
Let $p:Z\to N$ be a Legendrian fibration of a co-oriented contact manifold.
Suppose that there exists a contact covering $ST^*M\to Z$ such that $M$ is open
and the pre-image of any fibre of $p$ is a union of fibres of $ST^*M$.
(Note that $p$ does not have to be locally trivial and its fibres do not have to be spheres.)
Then Theorem~\ref{mainresult} and Corollary~\ref{loops} hold for the fibres of $p$.
\end{rem}

\subsection{Orderability of $ST^*M$}
\label{orderapplic}
Let $(Y,\ker\alpha)$ be a connected contact manifold. Consider the identity
component ${\mathrm{Cont}}_0(Y)$  of the group of compactly
supported contactomorphisms of $(Y,\ker\alpha)$ and let
$\widetilde{\mathrm{Cont}}_0(Y)$ denote the universal cover
of this group corresponding to the base point ${\mathrm{id}_Y}\in {\mathrm{Cont}}_0(Y)$.
For $f,g\in\widetilde{\mathrm{Cont}}_0(Y)$, write $f\preceq g$ if
the element $gf^{-1}$ can be represented by a path $\phi_t\in {\mathrm{Cont}}_0(Y)$
such that the contact Hamiltonian $H:=\alpha(\frac{d\phi_t}{dt})$ is non-negative.
Following Eliashberg and Polterovich~\cite{ElPo}, we say that the contact
manifold $(Y,\ker\alpha)$ is {\it orderable\/} if the relation $\preceq$ defines
a genuine partial order on $\widetilde{\mathrm{Cont}}_0(Y)$.

Eliashberg, Kim, and Polterovich~\cite[Theorem 1.18]{ElKiPo} used contact homology to prove that $ST^*M$
is orderable for a closed manifold $M$ such that its fundamental group $\pi_1(M)$
is either finite or has infinitely many conjugacy classes. It is an open problem
whether an infinite finitely presented group can have finitely many conjugacy classes,
see~\cite[Problem (FP19)]{BMS}. The following result shows that the orderability of $ST^*M$
does not depend on the solution of that problem.

\begin{cor}
\label{order}
$ST^*M$ is orderable for any closed manifold $M$.
\end{cor}

\begin{proof}
By~\cite[Criterion 1.2.C]{ElPo}, a closed contact manifold $(Y,\ker\alpha)$
is orderable if and only if there does not exist a contractible loop of contactomorphisms
$\phi_t\in\mathrm{Cont}_0(Y)$, $t\in [0,1]$, such that $\phi_0=\phi_1=\mathrm{id}_Y$
and the corresponding contact Hamiltonian is everywhere positive. It is clear that applying
a contact isotopy of $\mathrm{id}_Y$ generated by a positive contact Hamiltonian
to any Legendrian submanifold $L\subset Y$, we obtain a positive Legendrian isotopy
of~$L$. Thus, if $Y$ is not orderable, then every Legendrian isotopy class contains
a positive (contractible) Legendrian loop.

Suppose that $ST^*M$ is not orderable. By \cite[Theorem 1.18]{ElKiPo},
the fundamental group of $M$ is infinite and hence the universal cover of $M$ is open.
In that case, however, the Legendrian isotopy class of the fibre of $ST^*M$
does not contain positive Legendrian loops by Corollary~\ref{loops},
a contradiction.
\end{proof}

\begin{exm}
The proof of Corollary~\ref{order} shows that if $\pi_1(M)$ is infinite,
then there are no positive loops in ${\mathrm{Cont}}_0(ST^*M)$, contractible or not.
On the other hand, the (co-)geodesic flow of the standard round metric
on the $m$-sphere $S^m$ defines a {\it non-contractible\/} positive loop
in ${\mathrm{Cont}}_0(ST^*S^m)$, cf.\ Example~\ref{Yxl}.\qed
\end{exm}

\begin{cor}
\label{cover}
Let $p:ST^*M\to Z$ be a contact covering of a closed contact manifold $Z$.
Then $Z$ is orderable.
\end{cor}

\begin{proof} Suppose that $Z$ is not orderable and argue by contradiction.
By~\cite[Criterion 1.2.C]{ElPo}, there exists a contractible loop of contactomorphisms
of $Z$ based at $\id_Z$ and generated by a positive contact Hamiltonian.
Since this loop is contractible, it lifts to a loop of contactomorphisms of $ST^*M$
with the same properties. If $M$ is closed, we conclude that $ST^*M$ is not orderable
by~\cite[Criterion 1.2.C]{ElPo}, which contradicts Corollary~\ref{order}.
If $M$ is open, then the argument from the proof of that corollary
shows that there exists a positive Legendrian loop based at a fibre
of $ST^*M$, which contradicts Corollary~\ref{loops}.
\end{proof}

\subsection{Legendrian Low Conjecture}
\label{lowapplic}
Here we give a very brief exposition of the relevant material from Lorentz geometry.
Further details and references may be found in~\cite{ChNe}.

Let $(\ss,g)$ be a globally hyperbolic spacetime with a smooth spacelike Cauchy surface $M\subset\ss$.
That is to say, $(\ss,g)$ is a connected time-oriented Lorentz manifold and $M$ is a smooth spacelike
hypersurface in $\ss$ such that every inextensible future directed curve in $\ss$ meets $M$ exactly once.
(The time orientation on $X$ is a continuous choice of the future and past hemicones $C^\uparrow_x$
and $C^\downarrow_x$ in the non-spacelike cone at each $x\in X$. A piecewise smooth curve $\gamma=\gamma(t)$
is called future directed if $\dot\gamma(t)\in C^\uparrow_{\gamma(t)}$ for all~$t$.)

\begin{df}
Two points $x,y\in X$ are called {\it causally related\/} if they can be connected
by a future or past directed curve.
\end{df}

Let $\mathfrak N$ be the set of all future directed non-parameterised null geodesics in $(\ss,g)$
or, in other words, the set of all light rays of our spacetime. $\mathfrak N$ has a canonical
structure of a contact manifold, see~\cite[pp.\,252--253]{NatarioTod}.
There is a contactomorphism
$$
\rho_M:\mathfrak N\overset{\simeq}{\longrightarrow} ST^*M
$$
that associates to a null geodesic $\gamma\in\mathfrak N$ the equivalence class of
the (non-zero) linear form $v\mapsto g(\dot\gamma,v)$ on $T_{\gamma\cap M}M$,
where $\dot\gamma$ is a future pointing tangent vector to $\gamma$ at~$\gamma\cap M$.

The set $\mathfrak S_x$ of all null geodesics passing through a point $x\in\ss$ is
a Legendrian sphere in $\mathfrak N$ called the {\it sky\/} of that point.
Note that two skies intersect if and only if the corresponding points lie on the same
null geodesic. Note also that $\rho_M(\mathfrak S_x)=ST^*_xM$ for any $x\in M$.

Since $X$ is connected, the skies of any two points are Legendrian isotopic in $\mathfrak N$.
However, Legendrian links formed by unions of disjoint skies may be quite
different.

A basic observation is that all Legendrian links $\mathfrak S_x\sqcup\mathfrak S_y$
corresponding to causally {\it un\/}related points $x,y\in\ss$ belong to the
same Legendrian isotopy class, see~\cite[Lemma~4.3]{ChNe}.
Let us denote this isotopy class of Legendrian links by $\mathcal U$
(as in {\it u\/}nrelated and {\it u\/}nlinked). A natural way to represent $\mathcal U$
is to pick the points $x$ and $y$ on the Cauchy surface $M$ so that $\rho_M$ identifies
$\mathfrak S_x\sqcup\mathfrak S_y$ with $ST^*_xM\sqcup ST^*_yM\subset ST^*M$.

Two skies $\mathfrak S_x,\mathfrak S_y\subset\mathfrak N$ are said to be {\it Legendrian linked\/}
if either $\mathfrak S_x\cap\mathfrak S_y\ne\varnothing$ or the Legendrian link
$\mathfrak S_x\sqcup\mathfrak S_y$ does {\it not\/} belong to~$\mathcal U$.
We have just seen that if $\mathfrak S_x$ and $\mathfrak S_y$ are Legendrian linked,
then the points $x$ and $y$ are causally related.

\begin{df}
The {\it Legendrian Low Conjecture\/} holds for a globally hyperbolic spacetime
if two points in it are causally related if and {\it only if\/} their skies are Legendrian linked.
\end{df}

\begin{rem}
The general problem of describing causal relations in terms of linking in the space
of null geodesics originates from the work of Robert Low that was apparently inspired
by a question raised by Penrose, see e.\,g.~\cite{Low0} and~\cite{LowNullgeodesics}.
The Legendrian Low Conjecture was explicitly stated by Nat\'ario and Tod~\cite[Conjecture 6.4]{NatarioTod}
in the case when the Cauchy surface $M$ is diffeomorphic to an open subset of $\R^3$.
\end{rem}

It was shown in our paper \cite{ChNe} that the Legendrian Low Conjecture holds for
any globally hyperbolic spacetime such that its Cauchy surface has a cover diffeomorphic
to an open subset of $\R^m$, $m\ge 2$. Using Theorem~\ref{mainresult} instead of
\cite[Corollary 6.2]{ChNe}, we can now extend our result to a wider class of spacetimes.

\begin{thm}
\label{thma}
The {\it Legendrian Low Conjecture\/} holds for any globally hyperbolic spacetime
such that the universal cover of its Cauchy surface is not compact.
\end{thm}

\begin{proof}
Let $x,y\in\ss$ be two points such that their skies are disjoint and there exists
a future directed curve connecting $x$ to~$y$.
By \cite[Proposition 4.2]{ChNe}, there exists a non-negative Legendrian isotopy
connecting $\mathfrak S_y$ to $\mathfrak S_x$.
Suppose that the link $\mathfrak S_x\sqcup\mathfrak S_y$ belongs to $\mathcal U$.
Then the link $\rho_M(\mathfrak S_x\sqcup\mathfrak S_y)\subset ST^*M$ is Legendrian
isotopic to a link formed by a pair of fibres of $ST^*M$.
Since Legendrian isotopic links are ambiently contactomorphic,
we obtain a non-negative Legendrian isotopy connecting
two different fibres of $ST^*M$, which contradicts Theorem~\ref{mainresult}.
Thus, $\mathfrak S_x$ and $\mathfrak S_y$ are Legendrian linked.
\end{proof}

Combining this theorem with Perelman's proof of the Poincar\'e conjecture, we see that the
Legendrian Low Conjecture holds for any $(3+1)$-dimensional globally hyperbolic
spacetime such that the universal cover of its Cauchy surface is not diffeomorphic to $S^3$.
(This result was obtained in~\cite{ChNe} by a more involved argument using the
full strength of the geometrisation conjecture.) On the other hand,
if the Cauchy surface is a quotient of~$S^3$, then the Legendrian Low Conjecture
may fail because of the following general construction, cf.~\cite[Example 3]{ChernovRudyak}.

\begin{exm}
If $(M,\ov g)$ is a Riemann $Y_\ell^x$-manifold (see~Example~\ref{Yxl}), then the Legendrian Low Conjecture
does not hold for the globally hyperbolic spacetime $(M\times\R,\ov g\oplus -dt^2)$.
Indeed, null geodesics in this spacetime have the form $\gamma(s)=(\ov\gamma(s),s)$,
where $\ov\gamma$ is a $\ov g$-geodesic on $M$ and $s$ is the natural parameter on~$\ov\gamma$.
In particular, $\rho_{M\times\{0\}}(\mathfrak S_{(x,\ell)})=ST^*_xM$ by the definition of
a $Y_\ell^x$-manifold. Thus, the skies $\mathfrak S_{(x',0)}$ and $\mathfrak S_{(x,\ell)}$
are {\it not\/} Legendrian linked if $x'\ne x\in M$. However, the points $(x',0)$ and $(x,\ell)$
in $M\times\R$ are causally related if $x'$ is sufficiently close to $x$ in $M$.\qed
\end{exm}

\begin{rem}
One can use Remark~\ref{importantremark} and the proof of \cite[Theorem C]{ChNe}
to show that if the universal cover of the Cauchy surface of a globally hyperbolic
spacetime $\ss$ is non-compact, then
the Legendrian links $\mathfrak S_x\sqcup\mathfrak S_y$ and $\mathfrak S_y\sqcup\mathfrak S_x$
are different for any pair of causally related points $x,y\in\ss$ with disjoint skies.
\end{rem}


\begin{thebibliography}{99}
\bibitem{BMS}
G.~Baumslag, A.~G.~Myasnikov, V.~Shpilrain, {\it Open problems in combinatorial group
theory. Second edition}, Combinatorial and geometric group theory (New York,
2000/Hoboken, NJ, 2001), Contemp. Math. {\bf 296},
Amer. Math. Soc., Providence, RI (2002), 1--38.
\bibitem{BerardBergery}
L.~B\'erard-Bergery,
{\em Quelques exemples de vari\'et\'es riemanniennes o\`u toutes les
g\'eod\'esiques issues d'un point sont ferm\'ees et de m\^ eme longueur,
suivis de quelques r\'esultats sur leur topologie},
Ann. Inst. Fourier (Grenoble) {\bf 27} (1977), 231--249.
\bibitem{Besse}
A.~L.~Besse, {\em Manifolds all of whose geodesics are closed},
with appendices by D.~B.~A.~Epstein, J.-P.~Bourguignon,
L.~B\'erard-Bergery, M.~Berger and J.~L.~Kazdan.
Ergebnisse der Mathematik und ihrer Grenzgebiete, 93. Springer-Verlag, Berlin-New York, 1978.
\bibitem{Bh} M. Bhupal, {\it A partial order on the group of contactomorphisms of
$\mathbb R\sp {2n+1}$ via generating functions}, Turkish J. Math. {\bf 25} (2001), 125--135.
\bibitem{Ch} Yu. V. Chekanov, {\it Critical points of quasifunctions, and generating families of Legendrian manifolds},
Funktsional. Anal. i Prilozhen.  {\bf 30}:2 (1996), 56--69 (Russian);
English transl. in Funct. Anal. Appl. {\bf 30}:2 (1996), 118--128.
\bibitem{ChNe}
V. Chernov, S. Nemirovski,
{\it Legendrian links, causality, and the Low conjecture},
Preprint {\tt arXiv:0810.5091}.
\bibitem{ChernovRudyak}
V.~Chernov (Tchernov), Yu.~Rudyak,
{\it Linking and causality in globally hyperbolic space-times},
Comm. Math. Phys. {\bf 279} (2008), 309--354.
\bibitem{CFP}
V. Colin, E. Ferrand, P. Pushkar,
{\it Positive isotopies of Legendrian submanifolds},
Preprint available at {\tt http://people.math.jussieu.fr/\verb'~'ferrand/publi/PIL.pdf}
\bibitem{ElGr}
Ya. Eliashberg, M. Gromov,
{\it Lagrangian intersection theory\/{\rm :} finite-dimensional approach},
Geometry of differential equations,  27--118,
Amer. Math. Soc. Transl. Ser. 2, {\bf 186}, AMS, Providence, RI, 1998.
\bibitem{ElKiPo}
Y. Eliashberg, S. S. Kim, L. Polterovich,
{\it Geometry of contact transformations and domains\/{\rm :}
orderability versus squeezing}, Geom. Topol. {\bf 10} (2006), 1635--1747;
{\it Erratum}, Geom. Topol. {\bf 13} (2009), 1175--1176.
\bibitem{ElPo}
Y. Eliashberg, L. Polterovich,
{\it Partially ordered groups and geometry of contact transformations},
Geom. Funct. Anal. {\bf 10} (2000), 1448--1476.
\bibitem{Go}
C. Godbillon,
{\it Feuilletages. \'Etudes g\'eom\'etriques},
Progress in Mathematics, 98. Birkh\"auser Verlag, Basel, 1991.
\bibitem{LaSi} F. Laudenbach, J.-C. Sikorav,
{\it Persistance d'intersection avec la section nulle au cours d'une isotopie hamiltonienne dans un fibr\'e cotangent},
Invent. Math. {\bf 82} (1985), 349--357.
\bibitem{Low0} R.~J.~Low,
{\it Causal relations and spaces of null geodesics}, PhD Thesis, Oxford University (1988).
\bibitem{LowNullgeodesics}
R.~J.~Low, {\em The space of null geodesics}, Proceedings of the
Third World Congress of Nonlinear Analysts, Part 5 (Catania, 2000).
Nonlinear Anal. {\bf 47} (2001), 3005--3017.
\bibitem{NatarioTod}
J.~Nat\'ario, P.~Tod,
{\em Linking, Legendrian linking and causality},
Proc.~London Math.~Soc.~(3) {\bf 88} (2004), 251--272.
\bibitem{Perelman1}
G.~Perelman, {\em The entropy formula for the Ricci flow and its
geometric applications}, Preprint {\tt math.DG/0211159}.
\bibitem{Perelman2}
G.~Perelman, {\em Ricci flow with surgery on three-manifolds},
Preprint {\tt math.DG/0303109}.
\bibitem{Perelman3}
G.~Perelman, {\it Finite extinction time for the solutions to the Ricci flow
on certain three-manifolds}, Preprint {\tt math.DG/0307245}.
\bibitem{Vi} C. Viterbo, {\it Symplectic topology as the geometry of generating functions},
Math. Ann. {\bf 292} (1992), 685--710.
\end{thebibliography}
\end{document}